\documentclass[12pt]{amsart}
\usepackage{amsmath,amscd}
\theoremstyle{plain}

\newtheorem{thm}{Theorem}
\newtheorem*{thm*}{Theorem A}
\newtheorem*{cor*}{Corollary B}
\newtheorem*{thm**}{Theorem C}
\newtheorem*{thm***}{Theorem D}

\newtheorem{lemma}[thm]{Lemma}

\newtheorem*{conjecture*}{Conjecture}
\newtheorem{proposition}[thm]{Proposition}

\newtheorem*{claim*}{Claim}

\numberwithin{equation}{section}
\numberwithin{thm}{section}

\newcommand{\tS}{\widetilde{S}}
\newcommand{\tp}{\tilde{p}}
\newcommand{\ts}{\tilde{s}}
\newcommand{\tso}{\widetilde{s_0}}
\newcommand{\so}{s_0}
\newcommand{\tc}{\tilde{c}}
\newcommand{\tgi}{\widetilde{g_i}}
\newcommand{\tphi}{\widetilde{\varphi_i}}
\newcommand{\Teich}{\mathcal{T}_S}
\newcommand{\Mod}{\mathfrak{M}_S}
\newcommand{\barTeich}{\bar{\mathcal{T}}_S}
\newcommand{\C}{\mathbb C}
\newcommand{\R}{\mathbb R}
\newcommand{\Z}{\mathbb Z}
\newcommand{\V}{\mathbb V}
\newcommand{\Id}{\mathbf{I}}
\newcommand{\Diff}{\mathsf{Diff}(S)}
\newcommand{\Diffo}{\mathsf{Diff}^0(S)}
\newcommand{\Diffp}{\mathsf{Diff}^+(S)}
\newcommand{\tDiff}{\mathsf{Diff}(\tS)}
\newcommand{\Meto}{\mathsf{Met}_{-1}(S)}
\newcommand{\Met}{\mathsf{Met}(S)}
\newcommand{\CC}{\mathfrak{C}}
\newcommand{\SLtC}{\mathsf{SL}(2,\C)}
\newcommand{\SLtr}{\mathsf{SL}(2,\R)}
\newcommand{\PSLtC}{\mathsf{PSL}(2,\C)}
\newcommand{\PSLtr}{\mathsf{PSL}(2,\R)}
\newcommand{\scc}{\mathcal{C}(S)}
\newcommand{\inn}{\mathsf{Inn}}
\newcommand{\SLthr}{\mathsf{SL}(3,\R)}
\newcommand{\rpt}{\R\mathsf{P}^2}

\newcommand{\real}{\mathop{\mathfrak{Re}}\nolimits}
\newcommand{\ext}{\mathop{\rm ext}\nolimits}
\newcommand{\Hom}{\mathop{\rm Hom}\nolimits}

\newcommand{\iso}{\mathop{\rm Iso}\nolimits}
\newcommand{\dist}{\mathop{\rm dist}\nolimits}
\newcommand{\ord}{\mathop{\rm ord}\nolimits}
\newcommand{\vol}{\mathop{\rm vol}\nolimits}
\newcommand{\ad}{\mathop{\rm ad}\nolimits}
\newcommand{\rk}{\mathop{\rm rk}\nolimits}
\newcommand{\tr}{\mathop{\rm Tr}\nolimits}
\newcommand{\dbar}{\bar\partial}
\newcommand{\lra}{\longrightarrow}
\newcommand{\Hopf}{\mathop{\rm Hopf}\nolimits}
\newcommand{\simrightarrow}{\buildrel \sim\over\longrightarrow}
\renewcommand{\thefootnote}{\fnsymbol{footnote}}

\newcommand{\mcg}{\pi_0\big(\Diff\big)}
\newcommand{\Out}{\mathsf{Out}\big(\pi_1(S)\big)}
\newcommand{\Aut}{\mathsf{Aut}\big(\pi_1(S)\big)}
\newcommand{\Inn}{\mathsf{Inn}\big(\pi_1(S)\big)}
\newcommand{\dt}{d_{\mathcal{T}}}
\newcommand{\kx}{\mathcal{K}(X)}
\newcommand{\kT}{\mathcal{K}(\Teich)}
\newcommand{\qf}{\mathcal{QF}_S}
\newcommand{\Min}{\mathsf{Min}}
\newcommand{\rhthree}{\mathsf{H}^3}
\newcommand{\rhtwo}{\mathsf{H}^2}

\begin{document}
\title[Twisted harmonic maps]
{Energy of Twisted Harmonic Maps
of Riemann Surfaces
}

\author[Goldman]{William M. Goldman}
\address{Department of Mathematics \\
University of Maryland \\
College Park, MD 20742}
\thanks{Goldman supported in part by NSF grants DMS-0103889 and
DMS0405605}
\email{wmg@math.umd.edu}

\author[Wentworth]{Richard A. Wentworth}
\address{Department of Mathematics \\
Johns Hopkins University \\
Baltimore, MD 21218}
\thanks{Wentworth supported in part by NSF grant DMS-0204496}
\email{wentworth@jhu.edu}

\date{\today}

\begin{abstract}
The energy of harmonic sections of flat bundles of nonpositively
curved (NPC) length spaces over a Riemann surface $S$ is a function
$E_\rho$ on Teichm\"uller space $\Teich$ which is a qualitative
invariant of the holonomy representation $\rho$ of $\pi_1(S)$.
Adapting ideas of Sacks-Uhlenbeck, Schoen-Yau and Tromba, we show that
the energy function $E_\rho$ is proper for any convex cocompact
representation of the fundamental group. More generally, if $\rho$ is
a discrete embedding onto a normal subgroup of a convex cocompact
group $\Gamma$, then $E_\rho$ defines a proper function on the
quotient $\Teich/Q$ where $Q$ is the subgroup of the mapping class
group defined by $\Gamma/\rho(\pi_1(S))$.  When the image of $\rho$
contains parabolic elements, then $E_\rho$ is not proper.
Using the theory of geometric tameness developed by Thurston and
Bonahon~\cite{Bonahon},
we show that if
$\rho$ is a discrete embedding into $\SLtC$, then $E_\rho$ is proper
if and only if $\rho$ is quasi-Fuchsian. These results are used to
prove that the mapping class group acts properly on the subset of
convex cocompact representations.

\end{abstract}$\hbox{}$
\maketitle\thispagestyle{empty}
\tableofcontents
\section*{Introduction}
Let $S$ a closed orientable smooth surface with $\chi(S)<0$ and
$G$ a Lie group. This paper discusses
an analytic invariant of a representation $\pi_(S)\xrightarrow{\rho}G$,
and applies to the action of the mapping class group $\mcg$ of $S$
on the space of representations $\Hom(\pi_1(S),G)/G$.

We assume $G$ is a reductive real algebraic group with maximal compact
subgroup $K$ and symmetric space $X=G/K$.
Suppose that $\rho$ is {\em reductive,\/} that is, a
representation whose image is Zariski dense in a reductive subgroup of
$G$. Then according to Corlette~\cite{Corlette}, for every conformal
structure $\sigma$ on $S$, there is a $\rho$-equivariant harmonic
map
\begin{equation*}
\tS \xrightarrow{f_{\rho,\sigma}} X,
\end{equation*}
which is unique up to isometries of $X$.
(Such an equivariant harmonic map is called a {\em twisted harmonic map.\/})
In particular its {\em energy\/}
\begin{equation*}
E_\rho(\sigma)  \in \R
\end{equation*}
is well-defined. Letting $\sigma$ vary over {\em Teichm\"uller space\/}
$\Teich$ defines a function
\begin{equation*}
\Teich \xrightarrow{E_\rho} \R.
\end{equation*}
The starting point of our paper is the following result:
\begin{thm*}\label{thm:A}
Suppose that $\rho$ is convex cocompact.
Then $E_\rho$ is a proper function on $\Teich$.
\end{thm*}
Recall that a discrete subgroup $\Gamma\subset G$ is {\em convex cocompact\/}
if there exists a $\Gamma$-invariant closed geodesically convex subset
$N\subset X$ such that $N/\Gamma$ is compact.  A representation $\rho$
is {\em convex cocompact\/} if $\rho$ is an isomorphism of $\pi_1(S)$
onto a convex cocompact discrete subgroup of $G$.

>From this theorem follows the example which motivated this study.
Let $\CC$ be the subset of $\Hom(\pi_1(S),G)/G$ consisting
of equivalence classes of convex cocompact representations.

\begin{cor*}
$\mcg$ acts properly on $\CC$.
\end{cor*}
When $G=\PSLtC$, a convex cocompact representation is {\em quasi-Fuchsian,\/}
that is a discrete embedding whose action on $S2= \partial\rhthree$
is topologically conjugate to the action of a discrete subgroup of
$\PSLtr$.
The corollary is just the known fact that $\mcg$ acts properly on
the space $\qf$ of quasi-Fuchsian embeddings.
Bers's simultaneous uniformization theorem~\cite{Bers} provides
a $\mcg$-equivariant homeomorphism
\begin{equation*}
\qf \longrightarrow \Teich\times \barTeich.
\end{equation*}
Properness of the action of $\mcg$  on $\Teich$ implies properness on
on $\qf$.

The basic idea goes back to work of
Sacks-Uhlenbeck~\cite{SacksUhlenbeck} and
Schoen-Yau~\cite{SchoenYau}. When $\rho$ is a Fuchsian representation
(corresponding to a hyperbolic structure on $S$), Tromba~\cite{Tromba}
proved that $E_\rho$ is proper and has a unique critical point
(necessarily a minimum). When $\rho$ is a quasi-Fuchsian
$\PSLtC$-representation, $E_\rho$ is proper.
Uhlenbeck~\cite{Uhlenbeck} gave an explicit criterion for when
$E_\rho$ has a unique minimum.  Generally $E_\rho$ admits more than
one critical point, for quasi-Fuchsian $\rho$.  This follows from the
existence of quasi-Fuchsian hyperbolic 3-manifolds
containing arbitrarily many minimal surfaces, as constructed by Joel Hass
and Bill Thurston (unpublished).

However, as first shown by Kleiner and Leeb~\cite{KL} (see also
Quint~\cite{Quint}), convex cocompactness is highly restrictive,
only interesting when $G$ has $\R$-rank one. A more general condition
guaranteeing properness of the action of $\mcg$ is given by the
notion of {\em Anosov representations\/} introduced by
Labourie~\cite{Labourie}.

We generalize these results in two directions. First, we extend the
results on isometric actions of surface groups to {\em
isometric actions on non-positively curved metric spaces\/} as
developed by Korevaar-Schoen~\cite{KS1,KS2}. Second, following a
suggestion of Bruce Kleiner, we consider embeddings of surface groups
onto {\em normal subgroups\/} of a convex cocompact group $\Gamma$ of
isometries of an NPC space.  The quotient group
$Q=\Gamma/\rho(\pi_1(S))$ acts on $\Teich$.  Since $E_\rho$ is
$Q$-invariant, it induces a function $E'_\rho$ on $\Teich/Q$ and we
show:

\begin{thm**}\label{thm:C}
The mapping
\begin{equation*}
\Teich/Q \xrightarrow{E'_\rho} \R
\end{equation*}
is proper.
\end{thm**}
This generalization was motivated by hyperbolic 3-manifolds fibering
over the circle. The hyperbolic 3-manifold determines a representation
$\rho$ of the fundamental group of the fiber surface $S$.
Furthermore the monodromy of the fibration determines an automorphism
$\phi$ of $\pi_1(S)$ such that $\rho$ is conjugate to $\rho\circ \phi$.
According to Thurston~\cite{Thurston} (see also Otal~\cite{Otal}), $\phi$
is a {\em pseudo-Anosov\/} or {\em hyperbolic\/} mapping class, and
generates a proper $\Z$-action on $\Teich$. In particular every orbit is
an infinite discrete subset of $\Teich$.
Since $E_\rho$ is $\phi$-invariant and constant on each
infinite discrete orbit, $E_\rho$ is not proper.
Kleiner observed that $E_\rho$ induces a proper map on the cyclic
quotient $\Teich/\langle\phi\rangle$.

Properness of the energy function fails for surface group
representations containing ``accidental parabolics". Such representations
are discrete embeddings mapping some nontrivial simple loop $c$ to
a parabolic isometry. One can find a sequence of bounded energy
mappings for which the conformal structures $\sigma$ degenerate
as to shorten $c$, contradicting properness. Using the theory
of geometric tameness developed by Thurston and Bonahon~\cite{Bonahon}
(see the recent papers of Agol~\cite{Agol},
Calegari-Gabai~\cite{CalegariGabai} and Choi~\cite{Choi}),
we obtain a sharp converse to Theorem~A for discrete embeddings into $\PSLtC$:

\begin{thm***} \label{thm:qf}
Let $\rho: \pi_1(S)\to \PSLtC$ be a discrete embedding.  Then $E_\rho$
is proper if and only if $\rho$ is convex cocompact (that is,
quasi-Fuchsian).
\end{thm***}


\section*{Acknowledgments}
This paper grew out of conversations with many mathematicians over a
period of several years. In particular we thank Ian Agol,
Francis Bonahon, Dick Canary, David Dumas, Cliff Earle,
Joel Hass, Misha Kapovich, Bruce Kleiner, Fran\c cois Labourie,
John Loftin, Yair Minsky, Rick Schoen, Bill Thurston, Domingo Toledo,
Karen Uhlenbeck, Mike Wolf, and Scott Wolpert for helpful discussions.
We also thank the referee for several helpful suggestions.

\section*{Notation and Terminology}
If $X$ is a metric space, we denote the distance function by $d_X$.
If $c$ is a curve in a length space $X$, we denote its length by
$L_X(c)$. We denote by $[a]$ the equivalence class of $a$, in various
contexts. We denote the identity transformation by $\Id$.
We shall sometimes implicitly assume a fixed basepoint $s_0\in S$
in discussing the fundamental group $\pi_1(S)$
and the corresponding universal covering space $\tS\longrightarrow S$.

Although it is more customary to define the mapping class group
by {\em orientation-preserving\/} diffeomorphisms, for our purposes
it seems more natural to consider all diffeomorphisms. Orientation-reversing
mapping classes induce anti-holomorphic isometries of $\Teich$,
which are nonetheless appropriate in our setting.

\section{Flat bundles and harmonic maps}

Let $S$ be a closed oriented surface with $\chi(S)<0$, and let
$\pi_1(S)$ be its fundamental group.
Let $(X,d)$ be a complete nonpositively curved length space
(an \emph{NPC space}) with isometry group $G$.

Choose a universal covering space $\tS\longrightarrow S$ with group of
deck transformations $\pi_1(S)$.  An isometric action of $\pi_1(S)$ on
$X$ is a homomorphism $\pi_1(S)\xrightarrow{\rho} G$,  where
$G$ is the  isometry group of $X$.
Such a homomorphism defines a {\em flat $(G,X)$-bundle $X_\rho$\/}
over $S$, whose total space is the quotient $\tS\times X$ by the
(diagonal) $\pi$-action by deck transformations on $\tS$ and by $\rho$
on $X$.  Since every flat bundle over a simply connected space is
trivial, a section over the universal covering space $\tS$ is the
graph of a mapping $\tS\longrightarrow X$.  Sections of $X_\rho$
correspond to $\rho$-equivariant mappings
\begin{equation*}
\tS\xrightarrow{u}X.
\end{equation*}
Since $X$ is contractible (see, for example,
Bridson-Haefliger~\cite{BridsonHaefliger}),
sections always exist.

An important case (and the only one treated in this paper)
occurs when $\rho$ is a {\em discrete embedding\/} (otherwise known
as a discrete faithful representation). Then $\rho$ maps
$\pi_1(S)$ isomorphically onto a discrete subgroup $\Gamma\subset G$ and
determines a properly discontinouous free isometric action
of $\pi_1(S)$ on $X$. The quotient $X/\Gamma$ is a NPC space
locally isometric to $X$ with fundamental group $\Gamma \cong \pi_1(S)$.
Indeed, the representation $\rho$ defines a preferred isomorphism
of $\pi_1(S)$ with $\pi_1(X_\rho)$, that is, a preferred homotopy class
of homotopy equivalences $S \longrightarrow X/\Gamma$.
Sections of the flat $(G,X)$-bundle $X_\rho$
correspond to maps in this homotopy class.

For us, a {\em conformal structure\/} on $S$ will be an
{\em almost complex structure\/}   $\sigma$ on $S$, that is,
an automorphism of the tangent bundle $TS$ satisfying $\sigma2=-\Id$.
A Riemannian metric $g$ is {\em in the conformal class of $\sigma$\/}
if and only if
\begin{equation*}
g(\sigma v_1,\sigma v_2)= g(v_1,v_2)
\end{equation*}
for tangent vectors $v_1,v_2$.

Choose a conformal structure $\sigma$ on $S$.
Let $\tS\xrightarrow{f} X$
be a  continuously differentiable $\rho$-equivariant mapping.
Its differential defines a continuous section $df$
of the vector bundle
$T^\ast\widetilde S\otimes f^\ast TX$,
Choose a Riemannian metric
$g$ on $S$ in the conformal class of $\sigma$.
Denote by $\tilde{g}$ its pullback to $\tS$ and
$\widetilde{dA}$ the corresponding area form on $\tS$.
Let  $\Vert\cdot \Vert_{g,X}$ denotes the Hilbert-Schmidt norm with respect to
the metric on $\tS$ induced by $g$ and the metric on $X$.
Define {\em energy density\/} of $f$ with respect to $g$ on $\tS$ as
\begin{equation*}
\tilde{e}(f)=\Vert df\Vert_{g,X}^2 \;dA,
\end{equation*}

The energy density $\tilde{e}(f)$ is a $\pi_1(S)$-invariant exterior
2-form on $\tS$  and hence defines
an exterior 2-form, the {\em energy density\/} $e(f)$ on $S$.
The {\em energy\/} $E_{\rho,g}(f)$ is the integral
\begin{equation*}
E_{\rho,g}(f) = \int_S {e}(f).
\end{equation*}
Alternatively, $E_{\rho,g}(f)$ is the integral of $\tilde{e}(f)$
on $\tS$ over a
fundamental domain for the $\pi_1(S)$-action on $\tS$.  Since $S$ is
two-dimensional, $E_{\rho,g}(f)$
depends only on the
conformal structure $\sigma$, and we denote it $E_{\rho,\sigma}(f)$.


When the target $X$ is only a metric space,  define the
energy density via
\begin{equation*}
\tilde e(f)=\lim_{\varepsilon\downarrow 0}\int_{d_{\widetilde S}(x,y)=\varepsilon} \frac{ d^2_X(f(x), f(y))}{\varepsilon^2}\frac{ds(y)}{2\pi\varepsilon}\  .
\end{equation*}
(See Korevaar-Schoen \cite{KS1} or Jost \cite{Jost1}.)  For finite
energy maps the energy density $e(f)$ is a well-defined measure which is absolutely
continuous with respect to Lebesgue measure.  The Radon-Nikodym
derivative plays the role of $\Vert df\Vert^2$.  For more details, see
Korevaar-Schoen~\cite{KS1}. Finite energy maps always exist.
Furthermore, energy minimizing sequences of uniformly Lipschitz
equivariant mappings exist (\cite{KS1}, Theorem 2.6.4).  In addition
to providing a definition of energy minimizing maps to metric spaces,
their construction defines a Sobolev completion of the continuously
differentiable maps to Riemannian targets which does not appeal to an
isometric embedding of $X$ into euclidean space.

In many cases the infimum of the energy is realized.  In the context
of NPC targets, recall that a map $f$ is called {\em harmonic\/} if it
minimizes $E_{\rho,\sigma}$ among all $\rho$-equivariant maps of
finite energy.

The fundamental existence theorem for harmonic maps to nonpositively
curved Riemannian manifolds is due to Eells-Sampson~\cite{ES}.  In the
{\em twisted\/} (that is, equivariant)
setting there are various conditions on $\rho$ which
guarantee existence.  When $X$ is a symmetric space of noncompact
type, $\rho$ is said to be \emph{reductive} if its Zariski closure has
trivial unipotent radical.  Existence of a twisted harmonic map for reductive
$\rho$ was proven by Corlette~\cite{Corlette},
Donaldson~\cite{Donaldson}, Labourie~\cite{Labourie}
and Jost-Yau \cite{Jost-Yau}.  A geometric
notion of reductivity involving stabilizers of flat totally geodesic
subspaces was used in \cite{Labourie} (see also Jost \cite{Jost1}).
Korevaar and Schoen \cite{KS2} introduced the notion of a \emph{proper
action} (not to be confused with the more standard use of the term
\emph{proper} below) which is the condition that the sublevel sets of
the displacement function associated to a generating set of $\pi_1(S)$
are bounded.  This condition guarantees the existence of an energy
minimizer when $X$ is a general NPC space (see also \cite{Jost2}).

\section{Bounded geometry}
Let $\gamma\in G$. Its  translation length $\vert\gamma\vert$
is defined by:

\begin{equation} \label{E:translationlength}
\vert\gamma\vert:=\inf_{x\in X} d(x,\gamma x)\ .
\end{equation}

\begin{lemma} \label{L:minlength}
Let $\Gamma\subset G$ be a convex cocompact discrete subgroup.
Then $\exists\varepsilon_0>0$ such that
$|\gamma|\geq\varepsilon_0$ for all $\gamma\in\Gamma\setminus\{\Id\}$.
\end{lemma}

\begin{proof}
Suppose not.  Then
$\exists\gamma_i\in\Gamma$ such that
$|\gamma_i|\neq 0$ for all $i$, and $|\gamma_i|\to 0$.
Let $N$ be a closed convex $\Gamma$-invariant subset such that $N/\Gamma$
is compact.
Since $N$ is convex and $\Gamma$-invariant,  $\exists x_i\in N$ such that
\begin{equation*}
d(x_i, \gamma_i x_i)\lra 0\ .
\end{equation*}
Since $N/\Gamma$ is compact, $\exists \lambda_i\in\Gamma$ and $x\in N$
such that, after passing to a subsequence, $\lambda_ix_i\to x$.  Set
$\tilde\gamma_i=\lambda_i\gamma_i\lambda_i^{-1}$.  Then
\begin{equation*}
d(\lambda_i x_i, \tilde \gamma_i \lambda x_i)\lra 0\ .
\end{equation*}
Properness of the action of $\Gamma$  near $x\in N$ implies that for only
finitely many $i$ does $|\tilde\gamma_i|=|\gamma_i|$
This contradicts the assumption that
\begin{equation*}
0\neq |\gamma_i|\to 0.
\end{equation*}
\end{proof}

\begin{lemma} \label{L:paths}
Let $\varepsilon_0$ satisfy Lemma~\ref{L:minlength}. Let
$\gamma_1,\gamma_2\in\Gamma$ and $x, y\in X$. If
\begin{itemize}
\item $d(x,y)<\varepsilon_0/2$;
\item $d(\gamma_1 x, \gamma_2 y)<\varepsilon_0/2$,
\end{itemize}
then $\gamma_1=\gamma_2$.
\end{lemma}
\begin{proof}
\begin{align*}
|\gamma_2^{-1}\gamma_1|& \leq d(\gamma_2^{-1}\gamma_1 x, x) \\
&= d(\gamma_1 x, \gamma_2 x) \\
&\leq d(\gamma_1 x, \gamma_2 y) + d(\gamma_2 y, \gamma_2 x)        \\
&= d(\gamma_1 x, \gamma_2 y) + d(x, y) \; < \; \varepsilon_0\ .
\end{align*}
Now apply Lemma \ref{L:minlength}.
\end{proof}

\section{Existence of harmonic maps} \label{existence}

\begin{proposition} \label{P:existence}
Suppose that $\pi_1(S)\xrightarrow{\rho} G$ is convex cocompact.
Then there exists a $\rho$-equivariant harmonic map
$\tS\xrightarrow{u}X$.
\end{proposition}

We deduce this proposition as an immediate corollary of the following more
general proposition, which we state here for later applications.

\begin{proposition} \label{P:existence1}
Suppose that $\pi_1(S)\xrightarrow{\rho} G$ is an embedding
onto a normal subgroup of a convex cocompact subgroup
$\Gamma\subset\iso(X)$ such that $\rho(\pi_1(S))$ has trivial
centralizer in $\Gamma$.
Then there exists a $\rho$-equivariant harmonic map
$\tS\xrightarrow{u}X$.
\end{proposition}

\begin{proof}

For any NPC space $X$ and compact surface $S$, there exists an {\em
energy minimizing\/} sequence $u_i$ of uniformly Lipschitz
$\rho$-equivariant mappings $\tS\longrightarrow X$
(Korevaar-Schoen~\cite[Theorem 2.6.4]{KS1}).  Let $N\subset X$ be a
$\rho$-invariant convex set such that $N/\Gamma$ is compact.
Projection $X\longrightarrow N$ decreases distances, and therefore
decreases energy.  Thus we may assume that the image of $u_i$ lies in $N$.

Fix any point $\tso\in \tS$ with image $\ts\in S$.
Since $N/\Gamma$ is compact, after passing to a subsequence,
$\exists \gamma_i\in\Gamma$
such that $v_i(\tso)$  converges to a point in $N$, where
\begin{equation*}
v_i:=\rho(\gamma_i)\circ u_i.
\end{equation*}
The $v_i$ are uniformly Lipschitz and $v_i(\tso)$ converges.
The Arz\'ela-Ascoli theorem implies
that a subsequence of $v_i$
converges uniformly on compact subsets of $\tS$.
Choose $\varepsilon_0>0$ satisfying Lemma~\ref{L:minlength}.
For each compact $K\subset\tS$, there exists
$I>0$ so that

\begin{equation} \label{E:estimate1}
d(v_i(w), v_j(w)) <\varepsilon_0/2
\end{equation}
whenever $i,j\geq I$ and $w\in K$.

Each $v_i$ is equivariant with respect to
$\rho_i=\rho\circ \inn_{\gamma_i}$,
where $\inn_{\gamma_i}$
denotes the inner automorphism of $\pi_1(S)$ defined by $\gamma_i$.
Fix $i,j\geq I$, and set
$x=v_i(\tso)$ and  $y=v_j(\tso)$.

Choose a finite generating set $\Pi\subset\pi_1(S)$.
Applying \eqref{E:estimate1} to the finite set $K= \Pi\tso$,
\begin{equation*}
d(\rho_i(c)x, \rho_j(c)y)  = d(v_i(c \tso), v_j(c  \tso))  < \varepsilon_0/2\
\end{equation*}
whenever $c\in\Pi$. Since
\begin{equation*}
d(x,y)=d(v_i( \tso), v_j( \tso))<\varepsilon_0/2\ ,
\end{equation*}
Lemma \ref{L:paths} implies $\rho_i(c)=\rho_j(c)$ for all $c\in\Pi$.
As $\Pi$ generates $\pi_1(S)$ it follows $\rho_i=\rho_j$ if $i,j\geq
I$.  Since $\rho$ is injective and the centralizer of $\pi_1(S)$ in
$\Gamma$ is trivial, $c_i=c_j$ for all $i,j\geq I$.  Therefore $u_i$
itself converges locally uniformly to the desired minimizer.
\end{proof}

\section{The action of $\Diff$ on $\Teich$}\label{sec:actionTeich}
For later use, as well as a perspective on the theme of this paper,
we summarize in this section general facts on the action of the
diffeomorphism group on the space of metrics. A good general
reference for this material is Tromba's book~\cite{Tromba}.

Denote by $\Diff$ the group of smooth diffeomorphisms of $S$ with the
$C^\infty$ topology.
Let $\Diffo$ denote the identity component
of $\Diff$, that is, the group of all diffeomorphisms isotopic to the
identity. The {\em mapping class group\/} of $S$ is the quotient
\begin{equation*}
\mcg = \Diff/\Diffo.
\end{equation*}

The mapping class group relates to $\pi_1(S)$ as follows.
Let $s_0\in S$ be a fixed basepoint.
A diffeomorphism $\phi$ determines an automorphism of the fundamental
group $\pi_1(S,s_0)$ if $\phi(s_0)=s_0$.
Let $\phi\in\Diff$. Although $\phi$ may not fix $s_0$,
it is isotopic to  one which fixes $s_0$, which we call $\phi_1$.
This isotopy describes a path $q_1$ from $\phi(s_0)$ to $s_0$..
Suppose $\phi_2$ is another diffeomorphism isotopic to $\phi$ which fixes
$s_0$, with corresponding path $q_2$ from $\phi(s_0)$ to $s_0$.
Then the automorphisms of $\pi_1(S,s_0)$ induced by $\phi_1$ and $\phi_2$
differ by the inner automorphism $\inn_\gamma$
where $\gamma\in\pi_1(S,s_0)$ is the homotopy class of the based
loop $q_1\star (q_2)^{-1}$ in $S$. There results a homomorphism
\begin{equation*}
\pi_0(\Diff) \longrightarrow \Out
\end{equation*}
where
\begin{equation*}
\Out :=  \Aut/\Inn
\end{equation*}
is the quotient of $\Aut$ by its normal subgroup of
inner automorphisms.
\begin{thm}[Dehn-Nielsen]\label{eq:out}
The homomorphism
\begin{equation*}
\pi_0(\Diff) \longrightarrow \Out
\end{equation*}
is an isomorphism.
\end{thm}
We shall henceforth pass freely between these two approaches of the
mapping class group.  This was first proved by
Nielsen~\cite{Nielsen1927} and Dehn (unpublished).  For proof and
discussion, see Stillwell~\cite{Stillwell}) and
Farb-Margalit~\cite{FarbMargalit}).

Denote by $\Met$ the space of smooth Riemannian metrics on $S$ with
the $C^\infty$ topology.  For any smooth manifold $S$, the natural
action of $\Diff$ on $\Met$ is proper (Ebin~\cite{Ebin} and Palais
(unpublished) in general, and Earle-Eels~\cite{EarleEels} in dimension
$2$).  In particular its restriction to the subspace $\Meto$ of
metrics of curvature $-1$ is also proper.

Then $\Diffo$ acts properly on $\Meto$. The quotient,
comprising isotopy classes of hyperbolic structures on $S$,
identifies with the Teichm\"uller space of $S$
\begin{equation*}
\Meto/\Diffo  \longleftrightarrow \Teich \qquad \qquad
\end{equation*}
and inherits an action of the mapping class group
The properness of the action of $\Diff$ on $\Meto$ implies the following
basic fact:
\begin{thm}
\label{thm:proper}
$\mcg$ acts properly on $\Teich$.
\end{thm}
Closely related is the existence of a $\Diff$-invariant Riemannian
metric (in the Fr\'echet sense) on $\Met$. This
induces the $\mcg$-invariant
{\em Weil-Petersson metric\/} on $\Teich$.
This metric is incomplete, but complete metrics (for example
the Finslerian Teichm\"uller metric) exist which are $\mcg$-invariant.
For later applications, all we need is some $\mcg$-invariant
metric $\dt$ on $\Teich$. (For a survey of invariant metrics on $\Teich$
see Wolpert's paper~\cite{WolpertThisVolume}
in this volume.)

Theorem~\ref{thm:proper} is commonly attributed to Fricke.  The
customary proof uses a different set of ideas, more directly related
to representations of the fundamental group. We briefly digress to
sketch these ideas.

The uniformization theorem identifies $\Teich$ with a component of the
space of conjugacy classes of discrete embeddings
$\pi_1(S)\longrightarrow \SLtr$.  Such a representation is determined
up to conjugacy by its {\em character\/}
\begin{align*}
\pi_1(S)& \xrightarrow{\chi_\rho} \R  \\
 c & \longmapsto  \tr\rho(c).
\end{align*}
Geometrically $\chi_\rho$ corresponds to the {\em marked length
spectrum\/} $\ell_\rho$ which associates to a free homotopy class of
oriented loops in $S$ the length of the closed geodesic on
$\rhtwo/\rho(\pi_1(S))$ in that homotopy class.  Homotopy classes of
oriented loops in $S$ correspond to conjugacy classes in
$\pi_1(S)$). Denote this set of conjugacy classes by
$\widehat{\pi_1(S)}$. The key point is that the marked length spectrum
\begin{equation*}
\widehat{\pi_1(S)}\xrightarrow{\ell_\rho} \R_+
\end{equation*}
is finite-to-one (a proper map, where
$\widehat{\pi_1(S)}$ is discretely topologized).

Choose a $\mcg$-invariant metric $\dt$ on $\Teich$.  An {\em isometric\/}
action on a locally compact metric space is proper if and only if
some (and hence every) orbit is discrete.
Therefore it suffices to prove that every
$\mcg$-orbit is discrete. Suppose that $\phi_n\in\Aut$ is a sequence
of automorphisms and $\rho$ is a representation such that its images
$\rho\circ\phi_n$ converge to a representation $\rho_\infty$.
Let $\Pi\subset\pi_1(S)$
be a finite generating set and choose $C$ sufficiently large so that
\begin{equation*}
\ell_{\rho_\infty}(\gamma)\le C
\end{equation*}
for $\gamma\in \Pi$.
Then
\begin{equation*}
A :=  \{ \gamma\in\pi_1(S) \mid
\ell_{\rho_\infty}(\gamma)\le C \}.
\end{equation*}
is a finite union of conjugacy classes in $\pi_1(S)$ containing $\Pi$.
Let $\epsilon>0$. Then $\exists I$ such that
\begin{equation*}
\ell_{\rho\circ\phi_i}(\gamma)\le C + \epsilon
\end{equation*}
for $i \ge I$ and $\gamma\in A$.
Since
\begin{equation*}
\ell_{\rho\circ\phi_i}(\gamma)  =
\ell_{\rho}(\phi_i(\gamma)),
\end{equation*}
the set $A$ is invariant under all $\phi_i\circ (\phi_j)^{-1}$ for $i,j\ge I$.
>From this one can prove that the set of equivalence
classes $[\phi_i]\in\Out$ for $i\ge I$ is finite, so that
the sequence $[\rho\circ\phi_i]$ is finite, as desired.

For further details,
see Abikoff~\cite{Abikoff}, \S 2.2,
Farb-Margalit~\cite{FarbMargalit},
Harvey~\cite{Harvey}, \S 2.4.1,
Buser~\cite{Buser}, \S 6.5.6 (p.156),
Imayoshi-Tanigawa~\cite{IT},\S 6.3,
Nag~\cite{Nag},\S 2.7, and
Bers-Gardiner~\cite{BersGardiner}, Theorem II.

\section{Properness of the energy function} \label{S:proof}

We now prove that for  $\rho$ convex cocompact,
the function $E_\rho$ on $\Teich$ is proper.
With little extra effort, we prove a more general theorem
(suggested by Bruce Kleiner), concerning homomorphisms
\begin{equation*}
\pi_1(S) \xrightarrow{\rho} \Gamma \subset G
\end{equation*}
where $\Gamma$ is convex cocompact and $\rho(\pi_1(S))$ is a normal
subgroup $\Gamma_1 \lhd \Gamma$. Furthermore we assume that the
centralizer of $\Gamma_1$ in $\Gamma$ is trivial.
Let
\begin{equation*}
\Gamma\xrightarrow{\psi} \Aut
\end{equation*}
be the homomorphism induced
by the inclusion $\Gamma_1\hookrightarrow \Gamma$ and the isomorphism
$\pi_1(S)\xrightarrow{\rho} \Gamma_1$.
As $\Gamma_1$ has trivial centralizer,
$\psi$ is injective.
Thus $\psi$  induces a
monomorphism
\begin{equation*}
Q\hookrightarrow \Out
\end{equation*}
where $Q := \Gamma/\Gamma_1$.
Hence $Q$ acts on $\Teich$ via \eqref{eq:out}.
Furthermore $E_\rho$ is $Q$-invariant and hence induces a
map $\Teich/Q \xrightarrow{E'_\rho} \R$

\begin{proposition} The map $\Teich/Q \xrightarrow{E'_\rho} \R$ is proper.
\end{proposition}

Suppose that $[\sigma_i]\in \Teich$
is a sequence whose image in $\Teich/Q$ diverges.
Suppose further that
\begin{equation*}
E_\rho([\sigma_i])\leq B
\end{equation*}
for some constant $B> 0$, and  all $i=1,2,\dots$.

Our assumption that the images of $[\sigma_i]$ diverge in $\Teich/Q$
means the following.
Choose any invariant $\mcg$-invariant metric $\dt$ on $\Teich$.
We may assume, for each $\eta\in Q$, that
\begin{equation}\label{eq:divergent}
\dt(\psi(\eta)[\sigma_i],[\sigma_j])\geq 1
\end{equation}
for $i\neq j$.

By \cite{KS1} the $\rho$-equivariant harmonic maps
\begin{equation*}
(\tS,\tgi)\xrightarrow{u_i} X
\end{equation*}
have a uniform Lipschitz constant $K$ (depending on $B$),
where $\tgi$ denotes the
hyperbolic metric on $\tS$ associated to $\sigma_i$.  In particular, given
a closed curve $c$ in $S$, choose a lift $\tc\subset\tS$
running from $\tso$ to $[c]\tso$, where $[c]\in\pi_1(S;s_0)$ denotes the deck
transformation corresponding to $c$. Denote the length of $c$ with
respect to the metric $g_i$ on $S$ by $L_i(c)$. Then
\begin{align}\label{eq:lengthestimate}
\vert\rho([c])\vert & \leq d(u_i(\tso), \rho([c]) u_i(\tso)) \notag\\
&  =d (u_i(\tso), u_i(\rho([c]) \tso))\notag \\
& \leq L_X(u_i(\tilde c))\notag \\ & \leq K L_i(c)\ .
\end{align}

Suppose that $c\subset \Sigma$ is any closed essential curve. Since
$\rho$ is injective, the isometry $\rho(c)$ is nontrivial. Let
$\varepsilon_0>0$ satisfy Lemma~\ref{L:minlength}.
Then \eqref{eq:lengthestimate} implies
\begin{equation*}
\ell_c(\sigma_i)\geq\varepsilon_0/K
\end{equation*}
where $\ell_c(\sigma)$ denotes the
{\em geodesic length function of $c$ with respect to $\sigma$,\/} that
is, the length of the unique closed geodesic freely homotopic to $c$
in the hyperbolic metric corresponding to $\sigma$.

Mumford's compactness theorem~\cite{Mumford} implies that the
conformal structures $[\sigma_i]$ project to a compact subset of the
{\em Riemann moduli space\/} $\Teich/\mcg$.  Thus $[\varphi_i]\in\mcg$
and $[\sigma_\infty]\in\Teich$ exist such that, after passing to a
subsequence,
\begin{equation}\label{eq:converges}
[\varphi_i][\sigma_i]\longrightarrow [\sigma_\infty].
\end{equation}

As $\Diff$ acts properly on the set of Riemannian
metrics (\S\ref{sec:actionTeich}),
representatives $g_i\in\Meto$ and $\varphi_i\in\Diff$ exist with
$\varphi_i(g_i)\longrightarrow g_\infty$, where $g_\infty$ denotes the
hyperbolic metric associated to $\sigma_\infty$.
Choose a base point $\tso\in \tS$  with image $\so\in S$.
We may assume that $\varphi_i(\so)=\so$.
Let $\tphi\in\tDiff$ be the unique lift of $\varphi_i$ such that
$\tphi(\tso)=\tso$.

The map
\begin{equation*}
v_i: = u_i\circ\tphi^{-1}: \tS \longrightarrow X
\end{equation*}
is harmonic with respect to the metric $\phi_i(g_i)$ and equivariant
with respect to the homomorphism
\begin{equation*}
\rho\circ (\varphi_i^{-1}):\pi_1(S)\longrightarrow G.
\end{equation*}

The maps $v_i$ are
uniformly Lipschitz with respect to the metric $\widetilde{g_\infty}$ on $\tS$
induced from the metric $g_\infty$ on $S$.
In particular the family $\{v_i\}$ is equicontinuous.

Since $N/\Gamma$ is compact, $\exists\gamma_i\in \Gamma$
such that all $\gamma_i v_i(\tso)$ lie in a compact subset of $X$.

By the Arz\'ela-Ascoli theorem, a subsequence of
\begin{equation*}
w_i := \gamma_i\circ v_i
\end{equation*}
converges
uniformly on compact sets.
For $I$ sufficiently large,
\begin{equation} \label{E:estimate}
\sup_{z\in\tS}\,d_{X}(w_i(z), w_j(z)) \;<\;\varepsilon_0/2\
\end{equation}
for $i,j\geq I$.

Each $v_i = u_i\circ \tphi^{-1}$
is equivariant with respect to $\rho\circ(\varphi_i)_\ast^{-1}$ and is
harmonic with respect to $\varphi_i(g_i)$.
Thus each $w_i = \gamma_i\circ v_i$
is equivariant with respect to
\begin{equation*}
\rho_i :=
\rho\circ(\varphi_i)_\ast^{-1}\circ\psi(\gamma_i)
\end{equation*}
and also harmonic with respect to $\varphi_i(g_i)$ (since $\gamma_i$ is
an isometry).

Fix $i,j\geq I$, and let
$x=w_i(\tso)$ and  $y=w_j(\tso)$.  For every $c\in \pi_1(S)$,
\eqref{E:estimate} implies
\begin{align*}
d_X(\rho_i(c)x, \rho_j(c)y) & = d_X(w_i(c\tso), w_j(c\tso)) \\
& <\; \varepsilon_0/2\ .
\end{align*}
Since $d_X(x,y) <\;\varepsilon_0/2$,
Lemma \ref{L:paths} implies $\rho_i(c)=\rho_j(c)$ for all
$c\in\pi_1(S)$.  Since $\rho$ is injective,
\begin{equation*}
\psi(\gamma_i)\circ(\varphi_i)_\ast=
\psi(\gamma_j)\circ(\varphi_j)_\ast
\end{equation*}
Theorem~\ref{eq:out},
implies the natural homomorphism
\begin{equation*}
\mcg \cong \Out
\end{equation*}
is injective; thus
$\psi(\gamma_i)\circ\varphi_i$
is isotopic to $\psi(\gamma_j)\circ\varphi_j$ for all $i,j\geq I$.
Call this common mapping class $[\varphi]$.
Thus, for $i\ge I$,
\begin{equation}\label{eq:stabilizing}
\psi(\gamma_i)\circ(\varphi_i)_* = [\varphi]
\end{equation}
If $i,j\geq I$, then
\begin{align*}
\dt\big(\,
\psi(\gamma_i)^{-1}[\sigma_i] \, ,\, &
\psi(\gamma_j)^{-1}[\sigma_j]   \big)\\  & =\;
 \dt\big(\,   [\varphi]^{-1}(\varphi_i^{-1})_*[\sigma_i],\,
    [\varphi]^{-1}(\varphi_j^{-1})_*[\sigma_j] \big)
\\ & \quad =\,
\dt\big(\,   (\varphi_i^{-1})_*[\sigma_i]\, ,\,
           (\varphi_j^{-1})_*[\sigma_j]    \big)
\\ & \qquad \lra 0
\end{align*}
by \eqref{eq:converges}, contradicting \eqref{eq:divergent}.
Thus $E_\rho$ is proper, as claimed.

\section{Action of the mapping class group}
Corollary~B follows from the properness of the
action of $\mcg$ on $\Teich$ and a general fact on proper actions on
metric spaces.
Let $X$ be a metric space and let $\kx$ denote
the space of compact subsets of $X$, with the Hausdorff metric.

\begin{lemma}\label{lem:compact}
A group $\Gamma$ of homeomorphisms of $X$ acts properly on $X$ if and only
if $\Gamma$ acts properly on $\kx$.
\end{lemma}
\begin{proof}
The mapping
\begin{align*}
X & \xrightarrow{\iota} \kx \\
x & \longmapsto \{x\}
\end{align*}
is a proper isometric $\Gamma$-equivariant embedding.
If $\Gamma$ acts properly on $\kx$, then equivariance implies
that $\Gamma$ acts properly on $X$.

Conversely, suppose that $\Gamma$ acts properly on $X$. For
any compact subset $K\subset\kx$ of $\kx$, its union
\begin{equation*}
UK := \bigcup_{A\in K} A
\end{equation*}
is a compact subset of $X$. For $\gamma\in\Gamma$ the
condition
\begin{equation}\label{eq:kk}
\gamma(K)\cap K\;\neq\;\emptyset
\end{equation}
implies the condition
\begin{equation}\label{eq:cupkk}
\gamma(UK)\,\cap\, UK \;\neq \;\emptyset.
\end{equation}
To show that $\Gamma$ acts properly on $\kx$, let $K\subset\kx$ be a
compact subset.  Since $\Gamma$ acts properly on $X$, only finitely
many $\gamma\in\Gamma$ satisfy \eqref{eq:cupkk}, and hence only
finitely many $\gamma\in\Gamma$ satisfy \eqref{eq:kk}. Thus $\Gamma$
acts properly on $\kx$.
\end{proof}

We now prove Corollary~B.
Let $[\rho]\in\CC$. By Theorem~A, $E_\rho$ is a proper function
on $\Teich$, and assumes a minimum $m_0(E_\rho)$
Furthermore
\begin{equation*}
\Min(\rho) := \{ [\sigma]\in\Teich \mid E_\rho(\sigma) = m_0(E_\rho) \}
\end{equation*}
is a  compact subset of $\Teich$, and
\begin{equation*}
\CC  \xrightarrow{\Min} \kT
\end{equation*}
is a $\mcg$-equivariant continuous mapping.

\begin{proof}[Conclusion of Proof of Corollary~B\;]
Lemmas~\ref{thm:proper} and \ref{lem:compact} together
imply $\mcg$ acts properly on $\kT$.  By
equivariance, $\mcg$ acts properly on $\CC$.\end{proof}

\section{Accidental parabolics} \label{S:parabolic}

Now we illustrate with a well-known construction how properness of the
energy functional can fail if the action contains non-semisimple
isometries.  For simplicity, assume in this section that $X$ is a
simply connected nonpositively curved complete Riemannian manifold (a
{\em Cartan-Hadamard manifold\/}) and $G$ its group of isometries.

\begin{thm}\label{thm:parabolics}
Let $\pi_1(S)\xrightarrow{\rho} G$ be a homomorphism.
Assume that for some simple closed curve $c$ in $S$,
there is a complete geodesic
\begin{equation*}
\R \xrightarrow{\gamma} X
\end{equation*}
and constants $C,\delta>0$ such that
\begin{equation}\label{eq:gamma}
d_X\big(\gamma(t)\, ,\,\rho[c]\gamma(t)\big)
\leq Ce^{-\delta t}\ ,
\end{equation}
for all $t\geq 0$.
Then the energy functional $E_\rho$ is not proper.
\end{thm}
\begin{proof}

It suffices to construct a family $\sigma_t$, $0<t\leq 1$, of
conformal structures on $S$ such that the corresponding points in
$\Teich$ diverge as $t\to 0$, and a family $u_t$ of $\rho$-equivariant
maps $\widetilde S\to X$ such that $E_{\rho,\sigma_t}(u_t)$ is
uniformly bounded in $t$.

Fix an initial conformal structure $\sigma_1$ on $S$.  Let
$A_\varepsilon$ denote a tubular neighborhood of the geodesic
representative of $c$ with respect to the hyperbolic metric $g_1$
associated to $\sigma_1$.  We denote this geodesic also by $c$.  Let
$A_\varepsilon^\pm$ be the connected components of $A_\varepsilon-c$.

We furthermore choose $A_\varepsilon$ such that in the uniformization
of $(S, g_1)$, $\widetilde A_\varepsilon^\pm$ are
isometric to the strip
\begin{equation*}
{\mathbb R}\,\times\,\big[\varepsilon,\,\frac1{\ell_c(\sigma_1)}\, \big)\,,
\end{equation*}
where $\varepsilon_1$ is some positive number,
and $\ell_c(g_1)$ denotes the length of the geodesic.  This realizes
the isometry $[c]\in \pi_1(S)$ as the isometry
\begin{equation*}
(x,y)\mapsto (x+1,y).
\end{equation*}

Define the family $\sigma_t$ of conformal structures by the
{\em plumbing construction\/} discussed by Wolpert~\cite{Wolpert}.
The conformal structure on the complement
$S_\varepsilon=S-A_\varepsilon$ remains fixed  whereas the
conformal structure on $ A_\varepsilon^\pm$ is equivalent to the
annulus
\begin{equation*}
A_t^\pm :=
\R/\Z \times[\varepsilon, 1/\ell_c(\sigma_t))
\end{equation*}
where $\ell_c(\sigma_t)\to 0$ as $t\to 0$.

Next, let $\gamma$ be the geodesic satisfying \eqref{eq:gamma}.
and let $W(t)$ be the quantity on the left-hand-side of \eqref{eq:gamma}:
\begin{equation*}
W(t) := d_X\big(\gamma(t)\, ,\,\rho[c]\gamma(t)\big).
\end{equation*}
Geodesically connect points on the geodesic $\gamma$ to the points
on its image $\rho[c]\gamma$ along geodesics as follows.
Define
\begin{equation*}
\R\times [0,\infty)\xrightarrow{\alpha} X
\end{equation*}
so that $s\mapsto \alpha(s,t)$ is
the complete unit speed geodesic satisfying
\begin{align*}
\alpha(0,t) & =\gamma(t) \\
\alpha(W(t),t) & =\rho[c]\gamma(t).
\end{align*}
Writing
\begin{equation*}
L(t)=(1/\delta)\log(t/\varepsilon),
\end{equation*}
notice that $W(L(t))\leq C\varepsilon/t$.  Define
\begin{align*}
[0,1]\times [\varepsilon,\infty) &\xrightarrow{\beta} X\ \\
(s,t) & \longmapsto \alpha\big(W(L(t)) s,L(t)\big)\ .
\end{align*}
Since
\begin{equation*}
\Vert
\partial_t \beta(0,t)\Vert=\frac{\Vert \alpha^\prime(L)\Vert t}{\delta}
= \frac{t}{\delta}.
\end{equation*}
the nonpositive curvature of $X$ implies
\begin{equation}\label{eq:dt}
\Vert \partial_t \beta(s,t)\Vert\leq t/\delta
\end{equation}
for all $0\leq s\leq 1$. Also,
\begin{equation*}
\Vert (\partial_s \alpha)(W(L(t)) s,L(t))\Vert = 1
\end{equation*}
so
\begin{align}\label{eq:ds}
\Vert \partial_s \beta(s,t)\Vert
& = W(L(t)) \Vert (\partial_s \alpha)(W(L(t)) s,L(t))\Vert \notag \\
& \leq C\varepsilon/t
\end{align}
Extend $\beta$ to ${\mathbb R}\times [\varepsilon,\infty)$
equivariantly with respect to the $\mathbb Z$-action
on $\R\times[\varepsilon,\infty)$
and $\rho(c)$ on $X$.
The derivative bounds  \eqref{eq:ds} and \eqref{eq:dt}
imply that $\beta$ has finite energy as an equivariant map.

Choose a finite energy $\rho$-equivariant map
$(\widetilde S,\sigma_1)\xrightarrow{u} X$.
The energy of its restriction $u_1^c$ to $\widetilde
S_\varepsilon$ is finite as well.  By interpolating near the
boundary $\partial S_\varepsilon$ we may assume that $u_1$ restricted
to the connected components of $\partial \widetilde S_\varepsilon $
coincides with the geodesic $s\mapsto \alpha(s,0)$.  Then for each
$t$, $u_1^c$ extends to a map $u_t:\widetilde S\to X$ by requiring
\begin{equation*}
u_t\bigr|_{\widetilde A_t^\pm}=\beta\bigr|_{{\mathbb R}\times
[\varepsilon, 1/\ell_c(\sigma_t)]}\ .
\end{equation*}

Then $u_t$ is equivariant and has finite energy with respect to
$\sigma_t$, uniformly in $t$.  This completes the proof.
\end{proof}

\section{When $G=\PSLtC$}

For discrete embeddings in $G = \PSLtC$, R.\ Canary and Y.\ Minsky
have explained a partial converse to Theorem~A. Namely suppose that
$\rho$ is a discrete embedding of a closed surface group $\pi_1(S)$ into
$G$. Let $M := \rhthree/\rho(\pi_1(S))$ be the corresponding
hyperbolic 3-manifold. We show (Theorem~C) that unless $\rho$ is
quasi-Fuchsian, then $E_\rho$ is not proper.
Assume that $\rho$ is not quasi-Fuchsian. Further assume that
$\rho(\pi)$ contains no parabolics; otherwise by
Theorem~\ref{thm:parabolics}, $E_\rho$ is not proper.

Under these assumptions, the work of Thurston and Bonahon~\cite{Bonahon}
guarantees a sequence of pleated surfaces
\begin{equation*}
S \xrightarrow{f_n} M
\end{equation*}
which exhaust the ends of the hyperbolic 3-manifold $M3$.

The intrinsic geometry of each $f_n$ is that of a totally geodesic
surface in $\rhthree$ and therefore its energy (computed with respect to
the intrinsic hyperbolic metric) equals
\begin{equation*}
-2\pi\chi(S)=\mathrm{area}(S).
\end{equation*}
Let $\sigma_n$ be the conformal structure underlying this intrinsic metric;
then
\begin{equation*}
E_\rho(\sigma_n) \le -2\pi\chi(S)
\end{equation*}
is bounded.

However, the corresponding sequence $[\sigma_n]\in\Teich$
tends to $\infty$. It suffices to show that for some $c\in\pi_1(S)$,
the geodesic length $\ell_c(\sigma_n)$ is unbounded.
Choose a nontrivial element $c\in\pi_1(S)$.
Since each pleated surface $f_n$ is an isometric map, it suffices to show
that the closed geodesics $c_n$ on $f_n$
become arbitrarily long. Otherwise, $\exists C$ such that
\begin{equation*}
\ell_{f_n}(c_n) \le C.
\end{equation*}
Let $c$ denote geodesic in $M$ corresponding to $\rho(c)$.
Since the pleated surfaces $f_n$ tend to $\infty$,
\begin{equation*}
d(c, f_n) \longrightarrow \infty
\end{equation*}
and in particular the curves $c_n$ (each homotopic to $c$)
become arbitrarily long, as claimed.

Thus, the energy function $E_\rho$ for a discrete embedding
$\pi_1(S)\xrightarrow{\rho}\SLtC$ is proper if and only if
$\rho$ is quasi-Fuchsian.

\section{Speculation}

Deformation spaces of flat bundles over a surface $S$ are natural
geometric objects upon which the mapping class group of $S$ acts.
When $G$ is a compact group, then the action is ergodic
(Goldman~\cite{Erg} and Pickrell-Xia~\cite{PickrellXia}).
At the other extreme, uniformization identifies the Teichm\"uller space
$\Teich$ of $S$ with a connected component in the deformation space
of flat $\PSLtr$-bundles over $S$, and $\mcg$ acts properly on $\Teich$.
In general one expects the dynamics of $\mcg$ to intermediate between
these two extremes.

As mentioned earlier, convex cocompactness excludes all higher rank
examples which do not come from rank one.  However it may be possible
to replace geodesic convexity of the Riemannian structure by another
notion. All that is needed is a compact {\em core\/} $N/\Gamma$ of the
locally symmetric space $X/\Gamma$ in which all all harmonic mappings
$S\longrightarrow X/\Gamma$ take values.

For example, when $G=\SLthr$, the mapping class group $\mcg$
acts properly on the component of $\Hom(\pi,G)/G$ corresponding to
convex $\rpt$-structures (Goldman~\cite{ConvexRPT}). Recently using
his notion of {\em Anosov representations,\/} Labourie has
proved~\cite{Labourie} that for any split real form $G$, the action of
$\mcg$ on the {\em Hitchin-Teichm\"uller component\/} of
$\Hom(\pi,G)/G$ (see Hitchin~\cite{Hitchin}) is proper.

Labourie's definition is as follows. The unit tangent bundle
\begin{equation*}
US\xrightarrow{\Pi}S
\end{equation*}
induces a central extension of fundamental groups
\begin{equation*}
\Z \longrightarrow
\pi_1(US) \xrightarrow{\Pi_*} \pi
\end{equation*}
where the center $\Z$ of $\pi_1(US)$ corresponds to the fundamental
group of the fibers of $\Pi$. A representation $\rho:\pi\longrightarrow G$
and a linear representation of $G$ on a vector space $V$ defines
a flat vector bundle
\begin{equation*}
V_\rho  \longrightarrow US
\end{equation*}
with holonomy representation
$\rho\circ \Pi_*$. Let $\tilde{\xi}_t$ denote the lift of the vector field
on $US$ defining the geodesic flow to the total space $\V_\rho$.
Labourie defines an {\em Anosov structure\/} to be a continuous splitting
of the vector bundle
\begin{equation*}
V_\rho = V_+ \oplus V_0 \oplus V_-
\end{equation*}
so that vectors in $V_+$ (respectively in $V_-$) are exponentially
expanded (respectively contracted) under $\tilde{\xi}_t$.

Labourie proves~\cite{Labourie} that the mapping class group acts
properly on all such representations.
All known examples of open sets of representations upon which the
mapping class group acts properly satisify Labourie's condition.
The key point is reminiscent of the proof of properness in
\S\ref{sec:actionTeich}: from the representation he constructs a
class function $\pi_1(S)\xrightarrow{\ell_\rho} \R_+$ which is
bounded with respect to length function for (any) hyperbolic structure
on $S$ (or the word metric on $\pi_1(S)$.

In another direction, using ideas generalizing those of
Bowditch~\cite{Bowditch} Tan, Wong and Zhang~\cite{TWZ} have shown
that the action of $\mcg$ on representations satisfying the analogue
of {\em Bowditch's Q-conditions\/} is proper. This also generalizes
the properness of the action on the space of quasi-Fuchsian
representations.

\end{document}